\documentclass[a4paper,onecolumn]{revtex4}

	\usepackage{amsmath}
	\usepackage{amsfonts}
	\usepackage[ansinew]{inputenc}
	\usepackage[usenames,dvipsnames]{pstricks}
	\usepackage{subfigure}
	\usepackage{epsfig}
	\usepackage{pst-grad} 
	\usepackage{pst-plot} 
	\usepackage[hyperindex]{hyperref}
	\usepackage{cleveref}
	\usepackage{scrextend}
	\usepackage{amsthm}
	\usepackage{soul}

\changefontsizes[15pt]{11pt}

	\newtheorem{theorem}{Theorem}
	\newtheorem{corollary}{Corollary}

\allowdisplaybreaks

\begin{document}

\title{Quadratic Transformations of Hypergeometric Function and Series with Harmonic Numbers}

\author{Martin Nicholson}

\begin{abstract}
  In this brief note, we show how to apply Kummer's and other quadratic transformation formulas for Gauss' and generalized hypergeometric functions in order to obtain transformation and summation formulas for series with harmonic numbers that contain one or two continuous parameters. We also give a generating function of the sequence $\frac{(a)_n (1-a)_n}{(n!)^2}H_n$ as a combination of Gauss hypergeometric function and elementary functions.
\end{abstract}

\maketitle
There is an extensive literature on application of Newton-Andrews method for finding identities with harmonic numbers
$$
H_n=\sum_{k=1}^n\frac{1}{k}
$$
and generalized Harmonic numbers 
$$
H_n^{(r)}=\sum_{k=1}^n\frac{1}{k^r}.
$$
For a general description of this method and partial overview of the literature one can consult the papers [\onlinecite{chu,choi,wang}].
The aim of this paper is to study transformation formulas that are obtained from quadratic transformations of hypergeometric function by integration using Euler's integral representation for Harmonic numbers
\begin{equation}\label{euler}
H_n=\int_0^1\frac{1-z^n}{1-z}dz.
\end{equation} 
We will use the standard notation for hypergeometric function
$$
{}_rF_{s}\left({a_1,\ldots a_r\atop b_1,\ldots,b_s};z\right)=\sum_{n=0}^\infty\frac{(a_1)_n\ldots (a_r)_n}{n!(b_1)_n\ldots (b_s)_n}z^n
$$
where $(a)_n=a(a+1)\ldots (a+n-1)$ is the Pochhammer symbol.

\begin{theorem}\label{1}
Let $a$ and $b$ be arbitrary complex numbers such that $a+b-1/2$ is not a negative integer. Then
$$
2\sum_{n=1}^\infty\frac{(2a)_n(2b)_n}{n!\left(a+b+1/2\right)_n}\frac{H_n}{2^n}=\sum_{n=1}^\infty\frac{(a)_n(b)_n}{n!\left(a+b+1/2\right)_n}H_n,
$$
$$
4\sum_{n=1}^\infty\frac{(2a)_n(2b)_n}{n!\left(a+b+1/2\right)_n}\frac{H_n^2+H_n^{(2)}}{2^n}=\sum_{n=1}^\infty\frac{(a)_n(b)_n}{n!\left(a+b+1/2\right)_n}(H_n^2+H_n^{(2)}).
$$
\end{theorem}

\noindent
{\it{Proof.}} Starting from the quadratic transformation formula 2.11.2 in [\onlinecite{erd}]
\begin{equation}\label{kummer}
{}_2F_1\left({2a,2b\atop a+b+\frac{1}{2}};z\right)={}_2F_1\left({a,b\atop a+b+\frac{1}{2}};4z(1-z)\right),\quad\text{Re}~z<\frac12,
\end{equation}
replace $z$ in equation \eqref{kummer} with $\frac{z}{2}$ and multiply both sides by $(1-z)^{c-1}$ to get 
$$
\sum_{n=1}^\infty\frac{(2a)_n(2b)_n}{n!\left(a+b+1/2\right)_n2^n}{(1-z)^{c-1}z^n}{}=\sum_{n=1}^\infty\frac{(a)_n(b)_n}{n!\left(a+b+1/2\right)_n}(1-z)^{c-1}\left(2z-z^2\right)^n,
$$
where now $z\in[0,1]$. This series converges uniformly. Next we integrate this series termwise with respect to $z$ from 0 to 1. The integral on the LHS is easy to calculate and equals $B(c,n+1)$. The integral on the RHS is
\begin{align*}
\int_0^1(1-z)^{c-1}\left(1-(2z-z^2)^n\right)dz&=\int_0^1(1-z)^{c-1}\left(1-(1-z)^2\right)^ndz\\
&=\int_0^1z^{c-1}\left(1-z^2\right)^ndz\\
&=\frac12 \int_0^1z^{\frac{c-1}{2}}\left(1-z\right)^ndz\\
&=\frac12 B\left(\tfrac{c+1}{2},n+1\right).
\end{align*}
After multiplication of both sides by $c$ we deduce the series transformation
$$
\sum_{n=1}^\infty\frac{(2a)_n(2b)_n}{n!\left(a+b+1/2\right)_n2^n}\frac{\Gamma(c+1)\Gamma(n+1)}{\Gamma(c+n+1)}=\sum_{n=1}^\infty\frac{(a)_n(b)_n}{n!\left(a+b+1/2\right)_n}\frac{\Gamma(\tfrac{c}{2}+1)\Gamma(n+1)}{\Gamma(\tfrac{c}{2}+n+1)}.
$$ 
Next differentiate both sides of this equation with respect to $c$ at $c=0$ using formulas
$$
\dfrac{d}{dc}\frac{\Gamma(c+1)\Gamma(n+1)}{\Gamma(c+n+1)}\bigg|_{c=0}=-H_n,
$$
$$
\dfrac{d^2}{dc^2}\frac{\Gamma(c+1)\Gamma(n+1)}{\Gamma(c+n+1)}\bigg|_{c=0}=H_n^2+H_n^{(2)},
$$
from which the results stated in the theorem follow immediately. \qed

\begin{corollary}\label{corol2} Let $\psi(x)=\frac{\Gamma'(x)}{\Gamma(x)}$ denote the digamma function. Then
$$
\sum_{n=1}^\infty \frac{(n+1)(2a)_n}{(a+\frac32)_n2^n}H_n=\left(a+\tfrac{1}{2}\right) \left(\psi(a+\tfrac{1}{2})-\psi (\tfrac{1}{2})\right).
$$
\end{corollary}
\noindent{\it{Proof.}} In the first formula of Theorem \ref{1} put $b=1$
\begin{equation}\label{11}
2\sum_{n=1}^\infty\frac{(2a)_n(n+1)}{\left(a+3/2\right)_n}\frac{H_n}{2^n}=\sum_{n=1}^\infty\frac{(a)_n}{\left(a+b+1/2\right)_n}H_n.
\end{equation}
The sum on the left equals
\begin{align}\label{12}
\nonumber\sum_{n=1}^\infty\frac{(a)_n}{\left(a+3/2\right)_n}H_n&=\dfrac{d}{dc}{}_2F_1\left({a,c\atop a+\tfrac32};1\right)\bigg|_{c=1}\\
\nonumber &=\dfrac{d}{dc}\frac{\Gamma(a+\tfrac32)\Gamma(\tfrac32-c)}{\Gamma(\tfrac32)\Gamma(a+\tfrac32-c)}\bigg|_{c=1}\\
&=2\left(a+\tfrac{1}{2}\right) \left(\psi(a+\tfrac{1}{2})-\psi (\tfrac{1}{2})\right).
\end{align}
After combining \eqref{11} and \eqref{12} the proof is complete.\qed

\begin{corollary}\label{corol1} With the same notations as in Corollary \ref{corol2} we have
$$
\sum_{n=1}^\infty\frac{(a)_n (1-a)_n}{(n!)^2}\frac{H_n}{2^n}=\frac{\sqrt{\pi}}{2\Gamma(1-\tfrac{a}{2})\Gamma(\frac{a+1}{2})}(\psi(1-\tfrac{a}{2})+\psi(\tfrac{a+1}{2})-\psi(1)-\psi(\tfrac12)).
$$
\end{corollary}
\noindent{\it{Proof.}} This formula may be deduced from Theorem \ref{1} as a corollary, but the easiest way to prove it is by differentiating the summation formula 2.8.51 in [\onlinecite{erd}]
$$
{}_2F_1\left({a,1-a\atop c+1};\frac12\right)=\frac{\Gamma \left(\frac{c}{2}+1\right) \Gamma \left(\frac{c+1}{2}\right)}{\Gamma \left(\frac{c-a}{2}+1\right) \Gamma \left(\frac{c+a+1}{2}\right)}
$$ 
with respect to $c$ at $c=0$. \qed

\noindent{\it{Examples.}}
\begin{equation}\label{ex1}
\sum_{n=1}^\infty\binom{2n}{n}^2\frac{H_n}{32^n}=\frac{\Gamma^2(\tfrac{1}{4})}{4\sqrt{\pi}}
\left(1-\frac{4\ln2}{\pi}\right)
\end{equation}
\begin{equation}\label{ex2}
\sum_{n=1}^\infty\frac{(3n)!}{(n!)^3}\frac{H_n}{54^n}=\frac{\Gamma^3(\tfrac{1}{3})}{ 2^{7/3}\pi }\left(\sqrt{3}-\frac{9 \ln3}{2 \pi }\right)
\end{equation}
\begin{equation}\label{ex4}
\sum_{n=1}^\infty\binom{2n}{n}^2\frac{H_{2n}}{32^n}=\frac{\Gamma^2(\tfrac{1}{4})}{8\sqrt{\pi}}
\left(1-\frac{3\ln2}{\pi}\right)
\end{equation}
\begin{equation}\label{ex5}
\sum_{n=1}^\infty\frac{(3n)!}{(n!)^3}\frac{H_{3n}}{54^n}=\frac{\Gamma^3(\tfrac{1}{3})}{ 2^{7/3}\pi }\left(\frac{1}{\sqrt{3}}+\frac{2 \ln2-3 \ln3}{2 \pi }\right)
\end{equation}
Series that contain both central binomial coefficients and Harmonic numbers are studied in [\onlinecite{kalm,boyadzhiev,tauraso,chu2}], and the series where additionally the central binomial coefficients are squared or cubed have been studied in [\onlinecite{coffey,moll,chen,liu,sofo,campbell,guillera,wan}]. 

Formulas \eqref{ex1},\eqref{ex2} are direct consequences of the Corollary \ref{corol1}. To derive \eqref{ex4} and \eqref{ex5} observe that differentiating Kummer's summation formula ([\onlinecite{erd}], 2.8.50)
$$
{}_2F_1\left({2a,2b\atop a+b+\frac{1}{2}};\frac12\right)=\frac{\Gamma \left(\frac{1}{2}\right) \Gamma \left(a+b+\frac{1}{2}\right)}{\Gamma \left(a+\frac{1}{2}\right) \Gamma \left(b+\frac{1}{2}\right)}
$$
with respect to $b$ yields the summation formula [\onlinecite{choi2}]
$$
\sum_{n=1}^\infty\frac{(2a)_n(2b)_n2^{-n}}{n!\left(a+b+\tfrac12\right)_n}\sum_{k=0}^{n-1}\left(\frac{2}{2b+k}-\frac{1}{a+b+\tfrac12+k}\right)=\frac{\Gamma \left(\frac{1}{2}\right) \Gamma \left(a+b+\frac{1}{2}\right)}{\Gamma \left(a+\frac{1}{2}\right) \Gamma \left(b+\frac{1}{2}\right)},
$$
from which by specializing $a=b=\tfrac14$ one obtains the sum
$$
\sum_{n=1}^\infty\binom{2n}{n}^2\frac{4 H_{2 n}-3 H_n}{32^n}=\frac{\Gamma^2(\tfrac{1}{4})}{4\sqrt{\pi}}
\left(\frac{6 \ln2}{\pi }-1\right).
$$
Together with \eqref{ex1} this allows to calculate the sum \eqref{ex4}. Formula \eqref{ex5} is derived in analogous manner.

Note that the series
$$
\sum_{n=1}^\infty\binom{2n}{n}^2\frac{H_n}{16^n}k^{2n}=K(\sqrt{1-k^2})+\frac{1}{\pi}K(k)\log\frac{k^2}{16(1-k^2)},
$$
$$
\sum_{n=1}^\infty\binom{2n}{n}^2\frac{H_{2n}}{16^n}k^{2n}=\frac12K(\sqrt{1-k^2})+\frac{1}{\pi}K(k)\log\frac{k}{4(1-k^2)}
$$
have closed form in terms of logarithms and complete elliptic integrals of the first kind (see [\onlinecite{moll}] for a summation equivalent to a certain linear combination of the two formulas above).

\begin{theorem}\label{generalsummation} 
The following generalization of Corollary \ref{corol1} holds
\begin{align}\label{general}
\nonumber&\sum_{n=1}^\infty\frac{(a)_n (1-a)_n}{(n!)^2}H_nx^n=\frac{\pi }{2 \sin \pi  a} \, _2F_1\left({a,1-a\atop 1};1-x\right)\\
&+\frac12\left(\psi(1-\tfrac{a}{2})+\psi(\tfrac{a+1}{2})-\psi(1)-\psi(\tfrac12)-\frac{\pi }{\sin \pi  a}-\log \frac{1-x}{x}\right)\, _2F_1\left({a,1-a\atop 1};x\right).
\end{align}
\end{theorem}
\noindent
{\it{Proof.}} The hypergeometric function $\, _2F_1(a,b;c;x)$ satisfies the differential equation ([\onlinecite{erd}], section 2.1.1)
$$
x(1-x)u''+(c-(a+b+1)x)u'-abu=0.
$$
Let $v(x)=-\tfrac{d}{dc}\, _2F_1(a,1-a;c;x)\big|_{c=1}$. Then $v(x)$ satisfies the differential equation
\begin{equation}\label{nonhomogeneous}
x(1-x)v''+(1-2x)v'-a(1-a)v=f(x),
\end{equation}
where $f(x)=\tfrac{d}{dx}\, _2F_1(a,1-a;1;x)$. Note that $v(x)$ is the series on the left side of \eqref{general}.

The homogeneous equation \eqref{nonhomogeneous} with $f(x)=0$ has two linearly independent solutions $\, _2F_1(a,1-a;1;x)$ and $\, _2F_1(a,1-a;1;1-x)$. One can check by direct substitution (for example using computer algebra systems) that partial solution of the non-homogeneous equation \eqref{nonhomogeneous} is
$$
v_0(x)=-\frac12\, _2F_1(a,1-a;1;x)\log\frac{1-x}{x} .
$$
Thus the general solution of equation \eqref{nonhomogeneous} is
$$
v(x)=v_0(x)+A\, _2F_1(a,1-a;1;x)+B\, _2F_1(a,1-a;1;1-x).
$$
The coefficients $A$ and $B$ are determined from two conditions 
$$
v(0)=0
$$
$$
v(\tfrac12)=\frac12\left(\psi(1-\tfrac{a}{2})+\psi(\tfrac{a+1}{2})-\psi(1)-\psi(\tfrac12)\right)\, _2F_1\left({a,1-a\atop 1};\frac12\right),
$$
where the second condition is the Corollary \ref{corol1}. The first condition may be replaced by the weaker condition of the cancellation of the logarithmic singularity at $x\to +0$. Due to the asymptotics
$$
\, _2F_1(a,1-a;1;1-x)=\frac{\sin\pi a}{\pi}\,\log\frac{1}{x}+O(1),
$$
the first condition implies that $B=\frac{\pi}{2\sin\pi a}$. $A$ is now easily determined from the second condition.\qed

It is known that when $a=\frac{1}{r}$, $r=2,3,4,6$ the hypergeometric function $\, _2F_1(a,1-a;1;x)$ is an elliptic integral ([\onlinecite{berndt}], ch. 33). As a result both $\, _2F_1(\frac1r,1-\frac1r;1;x)$ and $\, _2F_1(\frac1r,1-\frac1r;1;1-x)$ have closed form in terms of gamma functions when $x$ are certain algebraic numbers, for example 
$$\, _2F_1\left(\tfrac{1}{3},\tfrac{2}{3};1;\tfrac{3(3-\sqrt{3})}{4} \right)=\sqrt{3} \, _2F_1\left(\tfrac{1}{3},\tfrac{2}{3};1;\tfrac{3 \sqrt{3}-5}{4}\right)=3^{3/8} (2+\sqrt3)^{1/4}\frac{\Gamma \left(\frac{1}{4}\right)^2}{(2 \pi )^{3/2}}.
$$
Thus Theorem \ref{generalsummation} allows one to generate closed form summation formulas with Harmonic numbers. Also note that for $r=2,3,4,6$ the expression $\frac{\left(1/r\right)_n \left(1-1/r\right)_n}{(n!)^2}$ takes the values $~\frac{1}{16^n}\binom{2n}{n}^2$,$~\frac{1}{27^n}\binom{3n}{2n}\binom{2n}{n}$,$~\frac{1}{64^n}\binom{4n}{2n}\binom{2n}{n}$,$~\frac{1}{432^n}\binom{6n}{3n}\binom{3n}{2n}$. It is also worth mentioning the summation formula
\begin{equation}
\sum_{n=1}^\infty\frac{(3n)!}{(n!)^3}\frac{H_{3n}}{27^n}x^n=\frac{\pi }{3 \sqrt{3}}\, _2F_1\left({\tfrac13,\tfrac23\atop 1};1-x\right)-\, _2F_1\left({\tfrac13,\tfrac23\atop 1};x\right)\log\frac{\sqrt{3 (1-x)}}{\sqrt[6]{x}}.
\end{equation}

Generating functions similar to \eqref{general} but with a sum of digamma functions instead of $H_n$ are considered in [\onlinecite{borwein}] as a special case of more general formulas from [\onlinecite{erd}], ch. 2.3.1.

\begin{theorem}\label{2}
Let $\mathrm{Re}(a+b)<\frac12$. Then
$$
\sum_{n=1}^\infty\frac{(2a)_n(2b)_n}{n!\left(a+b+1/2\right)_n}\frac{H_n}{n+1}=\frac{(2a+2b-1)\sin\pi a\sin\pi b}{(2a-1)(2b-1)\cos\pi(a+b)}\left(\psi(\tfrac12)+\psi(\tfrac32-a-b)-\psi(1-a)-\psi(1-b)\right).
$$
\end{theorem}
\noindent
{\it{Proof.}} Consider formula 2.11.7 in [\onlinecite{erd}]
$$
\frac{2\Gamma(\frac12)\Gamma(a+b+\frac12)}{\Gamma(a+\frac12)\Gamma(b+\frac12)}{}_2F_1\left({a,b\atop \frac{1}{2}};z^2\right)={}_2F_1\left({2a,2b\atop a+b+\frac{1}{2}};\frac{1+z}{2}\right)+{}_2F_1\left({2a,2b\atop a+b+\frac{1}{2}};\frac{1-z}{2}\right).
$$
We take the difference at $z=1$ and at $z$
$$
\frac{2\Gamma(\frac12)\Gamma(a+b+\frac12)}{\Gamma(a+\frac12)\Gamma(b+\frac12)}\sum_{n=1}^\infty\frac{(a)_n (b)_n}{n! \left(1/2\right)_n}(1-z^{2n})=\sum_{n=1}^\infty\frac{(2 a)_n (2 b)_n}{n! \left(a+b+1/2\right)_n}\left(1-\left(\frac{1-z}{2}\right)^n-\left(\frac{1+z}{2}\right)^n\right),
$$
then divide by $1-z^2$ and integrate. The integral on the LHS is
$$
\int_0^1\frac{1-z^{2n}}{1-z^2}dz=\sum_{k=1}^{n}\frac{1}{2k+1}=H_{2n}-\frac12 H_n,
$$
while the integral on the RHS
\begin{align*}
\int_0^1\left(1-\left(\frac{1-z}{2}\right)^n-\left(\frac{1+z}{2}\right)^n\right)\frac{dz}{1-z^2}&=\frac12\int_{-1}^1\left(1-\left(\frac{1-z}{2}\right)^{n-1}\right)\frac{dz}{1+z}\\
&=\frac12\int_0^1\frac{1-t^{n-1}}{1-t}dt=\frac12 H_{n-1}.
\end{align*}
Thus
$$
\sum_{n=1}^\infty\frac{(2 a)_n (2 b)_n}{n! \left(a+b+1/2\right)_n}H_{n-1}=\frac{2\Gamma \left(\frac{1}{2}\right) \Gamma \left(a+b+\frac{1}{2}\right)}{\Gamma \left(a+\frac{1}{2}\right) \Gamma \left(b+\frac{1}{2}\right)}\sum_{n=1}^\infty\frac{(a)_n (b)_n}{n! \left(1/2\right)_n}(2H_{2n}-H_n).
$$
To compute the sum on the RHS one can differentiate Gauss' summation formula ([\onlinecite{erd}], formula 2.8.46)
\begin{equation}\label{gauss}
{}_2F_1\left({a,b\atop c};1\right)=\frac{\Gamma(c)\Gamma(c-a-b)}{\Gamma(c-a)\Gamma(c-b)}
\end{equation}
with respect to $c$ at $c=\frac12$ to obtain
$$
-\sum_{n=1}^\infty\frac{(a)_n (b)_n}{n! \left(1/2\right)_n}(2H_{2n}-H_n)=\frac{\Gamma(\frac12)\Gamma(\frac12-a-b)}{\Gamma(\frac12-a)\Gamma(\frac12-b)}\left(\psi(\tfrac12)+\psi(\tfrac12-a-b)-\psi(\tfrac12-a)-\psi(\tfrac12-b)\right).
$$
After simplifying the product of Gamma functions using Euler's reflection formula, and after change of parameters $a\to a-\frac12$,$~b\to b-\frac12$ one recovers the formula stated in the theorem.\qed

Summation formulas of the same type as in Theorem \ref{2} were found in [\onlinecite{wei}] by applying Newton-Andrews method to analogs of Watson's sum
\begin{equation}\label{watson}
{}_3F_2\left({2a,2b,c \atop a+b+\frac12,2c};1\right)=\frac{\Gamma(\frac12)\Gamma(a+b+\frac12)\Gamma(c+\frac12)\Gamma(\frac12-a-b+c)}{\Gamma(a+\frac12)\Gamma(b+\frac12)\Gamma(\frac12-a+c)\Gamma(\frac12-b+c)}.
\end{equation}
Two of such analogs can be compactly written as [\onlinecite{lewanowicz}]
\begin{align}\label{watsonpm}
\nonumber \, _3F_2\left({2 a,2 b,c+\frac{\varepsilon}{2}\atop a+b+\frac{1}{2},2 c};1\right)&=\frac{\Gamma \left(\frac{1}{2}\right) \Gamma (c) \Gamma \left(a+b+\frac{1}{2}\right) \Gamma (c-a-b)}{\Gamma \left(a+\frac{1}{2}\right) \Gamma \left(b+\frac{1}{2}\right) \Gamma (c-a) \Gamma (c-b)}\\
&+\varepsilon \frac{\Gamma \left(\frac{1}{2}\right) \Gamma (c) \Gamma \left(a+b+\frac{1}{2}\right) \Gamma (-a-b+c)}{\Gamma (a) \Gamma (b) \Gamma \left(-a+c+\frac{1}{2}\right) \Gamma \left(-b+c+\frac{1}{2}\right)},\qquad \varepsilon=\pm 1.
\end{align}

Below we give a sketch of an alternative proof of Theorem \ref{2} using Newton-Andrews method. Differentiating \eqref{watson} with respect to $c$ at $c=1$ yields the sum
$$
\sum_{n=0}^\infty\frac{(2a)_n(2b)_n}{n!\left(a+b+1/2\right)_n}\left(\frac{2}{n+1}-\frac{2}{(n+1)^2}-\frac{H_n}{n+1}\right)
$$
in terms of gamma and digamma functions. The first term in the brackets is summed by \eqref{watson} with $c=1$. The sum with the second term can be expressed as
$$
\sum_{n=0}^\infty\frac{(2a)_n(2b)_n}{n!\left(a+b+1/2\right)_n}\frac{-2}{(n+1)^2}=-\frac{2a+2b-1}{(2a-1)(2b-1)}\lim_{c\to 0}\frac{1}{c}\left({}_3F_2\left({2a-1,2b-1,c \atop a+b-\frac12,2c+1};1\right)-1\right).
$$
Now the hypergeometric function is calculated using \eqref{watsonpm} with $\varepsilon=-1$, therefore it is possible to express this limit in terms of gamma functions and its derivatives. Hence the sum in Theorem \ref{2} is expressible in terms of gamma and digamma functions, as required.

\begin{theorem}\label{3}
Let $\mathrm{Re}\,(a+b)>0$, then
$$
\sum_{n=1}^\infty\frac{(1/2)_n(a+b)_n}{\left(1+a\right)_n(1+b)_n}H_{n}=4\sum_{n=1}^\infty\frac{(1-a)_n(1-b)_n}{(1+a)_n\left(1+b\right)_n}(-1)^n{H_n}+{}_3F_2\left({1,1/2,a+b\atop 1+a,1+b};1\right)\ln 4.
$$
\end{theorem}\noindent{\it{Proof.}} In the quadratic transformation formula 4.5.1 from [\onlinecite{erd}]
$$
{{}_{3}F_{2}}\left({a,b,c\atop a-b+1,a-c+1};-z\right)=(1+z)^{-a}{{}_{3}F_{2}}
\left({a-b-c+1,\frac{1}{2}a,\frac{1}{2}(a+1)\atop a-b+1,a-c+1};\frac{4z}{(1+z
)^{2}}\right)
$$
we put $a=1$ and then take the difference at $z=1$ and at $z$ to obtain
$$
\sum_{n=1}^\infty\frac{(b)_n(c)_n}{(2-b)_n\left(2-c\right)_n}(-1)^n(1-z^n)=\sum_{n=1}^\infty\frac{(1/2)_n(2-b-c)_n}{\left(2-b\right)_n(2-c)_n}\left(\frac12-\frac{1}{1+z}\left[\frac{4z}{(1+z
)^{2}}\right]^n\right).
$$
After dividing by $1-z$ and integrating termwise we get the following integral on the RHS
\begin{align*}
\int_0^1\left(\frac12-\frac{1}{1+z}\left[1-\left(\frac{1-z}{1+z}\right)^2\right]^n\right)\frac{dz}{1-z}&=\int_0^1\left(\frac12-\frac{1+t}{2}\left(1-t^2\right)^n\right)\frac{dt}{t(1+t)}\\
&=\frac12\int_0^1\left(\frac{1-\left(1-t^2\right)^n}{t}-\frac{1}{1+t}\right)dt\\
&=\frac14(H_n-\ln 4).
\end{align*}
As a result
$$
\sum_{n=1}^\infty\frac{(b)_n(c)_n}{(2-b)_n\left(2-c\right)_n}(-1)^nH_n=\frac14 \sum_{n=1}^\infty\frac{(1/2)_n(2-b-c)_n}{\left(2-b\right)_n(2-c)_n}(H_n-\ln4).
$$
One can easily bring this to the symmetric form stated in the theorem.\qed

By setting $a=\frac12,\, b=\frac14$ in theorem $4$ we get the following curious sum with harmonic numbers
\begin{corollary}\label{corol3}
\begin{equation}
\sum_{n=1}^\infty\frac{(\tfrac34)_n}{(\tfrac54)_n}\frac{\tfrac14-(-1)^n}{2n+1}H_n=\frac{\Gamma^4(\tfrac14)}{64\pi}\ln 2.
\end{equation}
\end{corollary}

\begin{theorem}\label{4}
Let $\mathrm{Re}\,b>\frac12$, then
$$
\frac14 \sum_{n=1}^\infty\frac{(1/2)_n(b)_n}{n!\left(2b\right)_n}{H_n}=\sum_{n=1}^\infty\frac{(1/2)_n(1-b)_n}{n!\left(b+1/2\right)_n}H_{2n}+\frac{\Gamma \left(b+\frac{1}{2}\right) \Gamma (2 b-1)}{\Gamma (b) \Gamma \left(2 b-\frac{1}{2}\right)}\ln 2 .
$$
\end{theorem}\noindent{\it{Proof.}} In the formula 2.11.5 from [\onlinecite{erd}]
$$
(1+z)^{-2a}{}_2F_1\left({a,b\atop 2b};\frac{4z}{(1+z)^2}\right)={}_2F_1\left({a,a+\frac12-b\atop b+\frac12};z^2\right)
$$
we put $a=\frac12$ and then consider its difference at $z=1$ and at $z$
$$
 \sum_{n=1}^\infty\frac{(1/2)_n(b)_n}{n!\left(2b\right)_n}\left(\frac12-\frac{1}{1+z}\left[\frac{4z}{(1+z
)^{2}}\right]^n\right)=\sum_{n=1}^\infty\frac{(1/2)_n(1-b)_n}{n!\left(b+1/2\right)_n}(1-z^{2n}).
$$
After dividing by $1-z$ and integrating with respect to $z$ one obtains $H_{2n}$ on the RHS, while the integral on the LHS was calculated in the proof of theorem \ref{3}. Thus
$$
\sum_{n=1}^\infty\frac{(1/2)_n(b)_n}{n!\left(2b\right)_n}H_n-\sum_{n=1}^\infty\frac{(1/2)_n(b)_n}{n!\left(2b\right)_n}\ln4=4\sum_{n=1}^\infty\frac{(1/2)_n(1-b)_n}{n!\left(b+1/2\right)_n}H_{2n}.
$$
Evaluating the second series using Gauss' sum \eqref{gauss} completes the proof.\qed

When $b$ is a positive integer Theorems \ref{3} and \ref{4} allow one to express the value of an infinite sum with Harmonic numbers as a finite sum.

{\it{Acknowledgements.}} The author of this paper wish to thank Dr. Wenchang Chu for valuable correspondence and encouragement.

\end{document}